\documentclass[11pt]{article}

\usepackage{enumerate}
\usepackage{amsmath,xspace,amssymb,mathrsfs}

\input xy
\xyoption{all}
\xyoption{2cell}
\UseAllTwocells
\CompileMatrices

\title{2-nerves for bicategories}
\author{Stephen Lack \thanks{The hospitality of Macquarie University and the 
support of the Australian Research Council are gratefully acknowledged.}\\ 
School of Computing and Mathematics\\ University of Western Sydney\\
email: {\tt s.lack@uws.edu.au} \and Simona Paoli \thanks{The support of the 
Australian Research Council is gratefully acknowledged.}\\ Department of 
Mathematics\\ Macquarie University\\ email:{\tt simonap@ics.mq.edu.au}}
\date{}

\renewcommand{\phi}{\varphi}
\renewcommand{\epsilon}{\varepsilon}

\newcommand{\A}{{\ensuremath{\mathscr A}}\xspace}
\newcommand{\B}{{\ensuremath{\mathscr B}}\xspace}
\newcommand{\C}{{\ensuremath{\mathscr C}}\xspace}
\newcommand{\E}{{\ensuremath{\mathbb E}}\xspace}
\newcommand{\G}{{\ensuremath{\mathscr G}}\xspace}
\renewcommand{\H}{{\ensuremath{\mathscr H}}\xspace}

\newcommand{\V}{{\ensuremath{\mathscr V}}\xspace}

\newcommand{\SSet}{\ensuremath{[\Delta\op,\Set]}\xspace}
\newcommand{\Set}{\textnormal{\bf Set}\xspace}
\newcommand{\cat}{\ensuremath{\textbf{Cat}_\textbf{1}}\xspace}
\newcommand{\Cat}{\textnormal{\bf Cat}\xspace}

\newcommand{\SCat}{\ensuremath{[\Delta\op,\Cat]}\xspace}
\newcommand{\CG}{\textnormal{\bf CG}\xspace}
\newcommand{\RCG}{\textnormal{\bf RCG}\xspace}

\newcommand{\twocat}{\textnormal{\bf 2-Cat}\xspace}
\newcommand{\Bicat}{\textnormal{\bf Bicat}\xspace}
\newcommand{\Bicats}{\textnormal{\ensuremath{\textbf{Bicat}_{\textbf{s}}}}%
\xspace}

\newcommand{\bicatl}{\textnormal{\ensuremath{\textbf{Bicat}_{\textbf{1}}}}%
\xspace}

\newcommand{\Deltac}{\ensuremath{\Delta_c}\xspace}
\newcommand{\Deltacop}{\ensuremath{\Delta\op_c}\xspace}
\newcommand{\Deltan}{\ensuremath{\Delta_n}\xspace}
\newcommand{\Deltanop}{\ensuremath{\Delta\op_n}\xspace}
\newcommand{\Deltatwo}{\ensuremath{\Delta_2}\xspace}
\newcommand{\Deltab}{\ensuremath{\Delta_b}\xspace}
\newcommand{\Deltatwoop}{\ensuremath{\Delta\op_2}\xspace}
\newcommand{\Deltathreeop}{\ensuremath{\Delta\op_3}\xspace}
\newcommand{\Deltabop}{\ensuremath{\Delta\op_b}\xspace}
\newcommand{\Jc}{\ensuremath{J_c}\xspace}

\newcommand{\Nc}{\ensuremath{N_c}\xspace}
\newcommand{\Ntwo}{\ensuremath{N_2}\xspace}
\newcommand{\Nb}{\ensuremath{N_b}\xspace}

\newcommand{\Cosk}{\textnormal{Cosk}\xspace}
\newcommand{\Ran}{\textnormal{Ran}\xspace}
\newcommand{\Lan}{\textnormal{Lan}\xspace}
\newcommand{\el}{\textnormal{el}\xspace}

\newcommand{\m}{\ensuremath{[m]}\xspace}
\newcommand{\n}{\ensuremath{[n]}\xspace}

\newcommand{\bd}[1]{\textbf{#1}\xspace}

\newcommand{\talg}{{\ensuremath{\textnormal{$T$-Alg}}}\xspace}
\newcommand{\talgs}{{\ensuremath{\textnormal{$T$-Alg}_{\text{s}}}}\xspace}
\newcommand{\Tam}{\textnormal{\bf Tam}\xspace}
\newcommand{\Simp}{\textnormal{\bf Simpson}\xspace}

\newcommand{\NHom}{\textnormal{\bf NHom}\xspace}
\newcommand{\NHoms}{\textnormal{\bf SHom}\xspace}
\newcommand{\Hom}{\textnormal{\bf Hom}\xspace}
\newcommand{\Catnps}{\ensuremath{\textnormal{\bf 2-Cat}_{\textnormal{\bf nps}}}\xspace}
\newcommand{\Catps}{\ensuremath{\textnormal{\bf 2-Cat}_{\textnormal{\bf ps}}}\xspace}
\newcommand{\Tamps}{\ensuremath{\textnormal{\bf Tam}_{\textnormal{\bf ps}}}\xspace}
\newcommand{\Tamnps}{\ensuremath{\textnormal{\bf Tam}_{\textnormal{\bf nps}}}\xspace}

\newcommand{\LxDbl}{\textnormal{\bf LxDbl}\xspace}

\newcommand{\Ps}{\textnormal{\bf Ps}\xspace}
\newcommand{\NPs}{\textnormal{\bf NPs}\xspace}

\newcommand{\two}{\ensuremath{{\hbox{\textrm 2}\kern-.25em
        \hbox{\vrule height1.5ex width 0.4pt depth -.2ex}}\kern.2em}\xspace}
\renewcommand{\t}{\times}
\newcommand{\op}{\ensuremath{^{\text{op}}}}

\newtheorem{theorem}{Theorem}[section]
\newtheorem{proposition}[theorem]{Proposition}   
   
\newtheorem{lemma}[theorem]{Lemma}   
\newtheorem{preremark}[theorem]{Remark}   
\newenvironment{remark}{\begin{preremark}\rm}{\end{preremark}}
\newtheorem{prexample}[theorem]{Example}

\newcommand{\proof}{\noindent{\sc Proof:}\xspace}
\def\endproof{~\hfill$\Box$\vskip 10pt}

\begin{document}

\label{firstpage}
\maketitle

\begin{abstract}
We describe a \Cat-valued nerve of bicategories, which associates to 
every bicategory a simplicial object in \Cat, called the 2-nerve. This
becomes the object part of a 2-functor $N:\NHom\to[\Delta\op,\Cat]$,
where \NHom is a 2-category whose objects are bicategories and whose 
1-cells are normal homomorphisms of bicategories.
The 2-functor $N$ is fully faithful and has a left biadjoint, and we
characterize its image. The 2-nerve of a bicategory is always a
weak 2-category in the sense of Tamsamani, and we show that \NHom is
biequivalent to a certain 2-category whose objects are Tamsamani weak
2-categories. 
\end{abstract}

This paper concerns a notion of ``2-nerve'', or \Cat-valued nerve, of 
bicategories.

To every category, one can associate its {\em nerve}; this is the simplicial 
set whose 0-simplices are the objects, whose 1-simplices are the morphisms, and
whose $n$-simplices are the composable $n$-tuples of morphisms. The face
maps encode the domains and codomains of morphisms, the composition law,
and the associativity property, while the degeneracies record information about
the identities.

This construction is the object part of a functor $N:\cat\to\SSet$ from the 
category of categories and functors, to the category of simplicial sets. This 
functor is fully faithful and has a left adjoint. It arises in a natural way, 
as the ``singular functor'' (see Section~\ref{sect:singular} below)
of the inclusion $J:\Delta\to\cat$
in \cat of the full subcategory $\Delta$ consisting of the non-empty finite 
ordinals. One can characterize the simplicial sets which lie in the image of 
the nerve functor as those for which certain diagrams are pullbacks. 

As observed by Street, one may define the nerve of a bicategory as the 
simplicial set whose 0-simplices are the objects, whose 1-simplices are the 
morphisms, whose 2-simplices consist of a composable pair $f$ and $g$ and a 
2-cell $gf\to h$, and so on. In this way, the category \bicatl of bicategories 
and {\em normal lax functors} becomes a full subcategory of \SSet; here a 
normal lax functor preserves the identities strictly, but preserves composition
only up to coherent, but not necessarily invertible, comparison maps. In the 
important special case of invertible comparison maps, one speaks rather of 
{\em normal homomorphisms}. Once again this nerve functor $\bicatl\to\SSet$ is 
a singular functor, this time of the inclusion $\Delta\to\bicatl$, where the 
non-empty finite ordinals are now seen as locally discrete (no non-identity 
2-cells) bicategories. The image of this nerve functor was
characterized explicitly in \cite{Duskin-nerve}.

In this construction, the 2-simplices are playing a double role: they encode 
both the 2-cells of the bicategory, and (indirectly) the composition of 
1-cells. It is essentially for this reason that the normal lax functors arise. 
In order to obtain a tighter control over the composition of 1-cells, and in 
particular to extract the normal homomorphisms, one could specify, as extra 
structure, which 2-simplices contain not just any 2-cell but an invertible one.
This gives rise to a structure called a {\em stratified} simplicial set,
which goes back to Roberts; see \cite{Verity-complicial} for a full
account, and a proof of the celebrated Street-Roberts conjecture, 
characterizing the nerves of $\omega$-categories.

Here we take a different approach. We consider only those 2-simplices in which 
the 2-cell is invertible; in order to compensate for this, however, we include 
arbitrary 2-cells as {\em morphisms} in a {\em category} of 1-simplices, and 
similarly there are categories of $n$-simplices for all higher $n$. The 
2-nerve of a bicategory is then a functor $X:\Delta\op\to\cat$, with $X_0$ the 
discrete (no non-identity 1-cells) category of objects of the bicategory; with 
$X_1$ the category whose objects are the morphisms of the bicategory, and 
whose morphisms are the 2-cells; and with $X_2$ a category whose objects 
consist of composable pairs $(f,g)$ with an {\em invertible} 2-cell $gf\to h$. 
The resulting 2-nerve construction is reminiscent of the {\em homotopy
coherent nerve} of \cite{Cordier-Porter}.

In order to describe this construction in terms of a singular functor, one 
needs to regard bicategories as objects of a {\em 2-category}. The morphisms
are the normal homomorphisms of bicategories, while the 2-cells are things
we call {\em icons}: these are oplax natural transformations 
$\alpha:F\to G:\A\to\B$ between normal homomorphisms, in which the component 
$FA\to GA$ is an identity, for each object $A$ of \A. We call the resulting
2-category \NHom. Once again, there
is a fully faithful inclusion $J:\Delta\to\NHom$, and the resulting nerve 
2-functor $N:\NHom\to\SCat$ is fully faithful (on both 1-cells and 2-cells). 
This $N$ does not have an adjoint in the ordinary strict sense, but it does 
have a left {\em biadjoint}, which is the most common situation for 
2-categories. Similarly \NHom does not have ordinary limits and colimits,
but it does have the bicategorical ones of \cite{FibBic}. These bicategorical 
results are obtained using the techniques of 2-dimensional universal algebra, 
as in \cite{BKP}. The image of $N$ can also be characterized, although this 
is somewhat more complicated than in the case of ordinary categories.

There are also variants \Hom, \Catps, and \Catnps, of \NHom in which one 
either generalizes from normal homomorphisms to homomorphisms, or specializes 
from bicategories to 2-categories, or both. All of these 2-categories are 
biequivalent. In particular, every bicategory is equivalent
in \NHom to a 2-category. (A normal homomorphism is an 
equivalence in \NHom if and only if it is bijective on objects and induces 
equivalences of hom-categories.)

We have chosen to take \NHom as basic, priveleging normal homomorphisms
over arbitrary homomorphisms. One reason for this is that it is straightforward
to ``normalize'' an arbitrary homomorphism so as to obtain a normal one;
indeed this arises in the fact, mentioned in the previous paragraph, that 
\NHom is biequivalent to \Hom. A second reason is that when we define our
2-nerves using $J:\Delta\to\NHom$, the resulting 2-nerves have a discrete
category of 0-simplices, which is part of Tamsamani's notion of weak 
2-category, mentioned below. If instead we had defined our 2-nerves using
the inclusion $\Delta\to\Hom$, the resulting categories of 0-simplices
would only have been equivalent to discrete categories, not, in general,
discrete. Our results could be adapted to that setting, but we find it
easier to restrict to the normal homomorphisms. Of course the composite
$\Hom\to\NHom\to[\Delta\op,\Cat]$ is still fully faithful in the bicategorical
sense: it induces equivalences (not isomorphisms) of hom-categories.


Tamsamani \cite{Tamsamani} and Simpson \cite{Simpson} have each defined 
notions of weak 2-category as
certain objects of \SCat; these notions then determine full sub-2-categories 
\Tam and \Simp, both containing the image of \NHom. The left 
biadjoint of the 2-nerve construction induces left biadjoints to the fully 
faithful inclusions $\NHom\to\Tam$ and $\NHom\to\Simp$, but these inclusions 
actually have 2-adjoints. The counits of these 2-adjunctions are invertible, 
while the units have components which are ``pointwise equivalences'': so for 
example if $X\in\Tam$, and $j:X\to NGX$ is the component at $X$ of the unit, 
then for each $\n\in\Delta$, the functor $j_n:X_n\to(NGX)_n$ is an equivalence 
of categories. A pointwise equivalence is not necessarily an equivalence, but 
it is always a ``weak equivalence''; this is enough to guarantee that if
one ``localizes the weak equivalences'', then \NHom, \Tam, and \Simp all 
become equivalent to the homotopy category of the Quillen model category 
\Bicat of \cite{qmbicat}. Rather than localizing, an alternative is to 
expand the notion of morphism of Tamsamani 2-categories: if one allows
not just 2-natural transformations, but pseudonatural ones, then the 
resulting 2-category \Tamps is in fact biequivalent to \NHom.

One of the motivations for this work was to determine the precise
relationship between bicategories and Tamsamani's weak 2-categories,
which was only very partially sketched in \cite{Tamsamani}. In fact our 
construction of the 2-nerve of a bicategory differs from that of 
\cite{Tamsamani}: see Remark~\ref{rmk:Tamsamani} below. On the other hand the 
bicategory we associate to a Tamsamani weak 2-category is the same as
the one constructed in \cite{Tamsamani}, however, unlike
\cite{Tamsamani}, we describe the functoriality of the
construction. This last point was also worked out in \cite{Paoli}, using a
slightly  different approach from that adopted here.

In Section~\ref{sect:singular} we describe two basic technical tools:
singular functors and coskeleta. In Section~\ref{sect:nerves} we
recall the basic facts about nerves of categories, while in 
Section~\ref{sect:2-nerves} we turn to 2-nerves. The next two sections
are not needed for the rest of the paper: the first describes our
basic 2-category \NHom of bicategories using 2-monads, and deduces
various useful things about it, while the second studies the
biequivalence between \NHom and various related 2-categories; in
particular, we see that every bicategory is equivalent in \NHom to a
2-category. In Section~\ref{sect:properties} we study various
properties of functors $\Delta\op\to\Cat$ which are 2-nerves of
bicategories, which leads to a characterization theorem in 
Section~\ref{sect:Tamsamani}, where we also establish the 2-adjunction
between \NHom and the 2-category of Tamsamani weak 2-categories, and
the precise relationship between these two structures. In particular,
we show that \NHom is biequivalent to \Tamps, the 2-category of 
Tamsamani weak 2-categories and pseudonatural morphisms.

The basic references for bicategories and 2-categories are still
\cite{bicategories,Kelly-Street,FibBic}.

\section{Singular functors and coskeleta}\label{sect:singular}

In this section we briefly recall some standard material on singular
functors and on coskeleta.
 
The results on singular functors are stated in terms of \V-categories,
for a symmetric monoidal closed \V which is complete and
cocomplete. The only cases needed will be the case 
$\V=\Set$ of ordinary categories, and the case $\V=\Cat$ of 2-categories.
With the exception of the second sentence of Proposition~\ref{prop:density},
everything here can be found in 
\cite[Chapter~5]{Kelly-book}. 

Let $F:\A\to\B$ be a functor with small domain. There is an induced functor 
$\B(F,1):\B\to[\A\op,\V]$ sending an object $B$ of \B to the functor 
$\B(F-,B):\A\op\to\V$, where $\B(F-,B)$ sends an object
$A$ to the hom-object $\B(FA,B)$. This functor $\B(F,1)$ is sometimes called 
the {\em singular functor} of $F$, and it may be obtained as the composite of 
the Yoneda embedding $\B\to[\B\op,\V]$ followed by the functor 
$[\B\op,\V]\to[\A\op,\V]$ given by restriction along $F$. 

The basic examples of such an $F$ will be the functor $J:\Delta\to\cat$
(where $\V=\Set$) and the 2-functor $J:\Delta\to\NHom$ (where $\V=\Cat$).
The resulting singular functors are then the nerve construction for 
categories and the 2-nerve construction for bicategories.

When $\B(F,1)$ is fully faithful, the functor $F$ is said to be dense. When 
$F$ is itself fully faithful, there is a characterization of when $F$ is dense,
involving \B being generated under colimits by \A, but we shall not need this 
characterization.
 
As observed by Kan, $\B(F,1)$ has a left adjoint provided that \B is 
cocomplete; the left adjoint can then be constructed as the left Kan extension 
of $F$ along the Yoneda embedding. It sends a presheaf $X:\A\op\to\V$ to the 
weighted colimit $X*F$, which may be given by the coend 
$$\int^A XA\cdot FA,$$
or, if $\V=\Set$, by the colimit of the functor
$$\el(X)\op \to \A \to \B.$$

We record for future reference the following, of which the first sentence
is \cite[Theorem~5.13]{Kelly-book}, and the second an easy consequence.

\begin{proposition}\label{prop:density}
If $F:\A\to\B$ and $G:\B\to\C$, with $G$ fully faithful, then $G$ is dense 
provided that $GF$ is so, and then the identity $GF=GF$ exhibits $G$ as
the left Kan extension of $GF$ along $F$. Furthermore, the singular 
functor $\C(G,1):\C\to[\B\op,\V]$ can then be obtained by first applying
the singular functor $\C(GF,1):\C\to[\A\op,\V]$ and then right Kan extending
along $F:\A\op\to\B\op$.
\end{proposition}

\proof
For the second sentence observe that $\Ran_{F}\C(GF,C)\cong\C(\Lan_{F}GF,C)$.
\endproof

The results on coskeleta are stated in terms of simplicial objects 
$[\Delta\op,\E]$ in a category \E with finite limits. Once again, the
main cases will be $\E=\Set$ and $\E=\cat$. Let $\Deltan$ be the full
subcategory of $\Delta$ consisting of all objects \m with $m\le n$, 
and $H_n:\Deltan\to\Delta$ the inclusion. The restriction along $H_n$
gives a functor $[\Delta\op,\E]\to[\Deltanop,\E]$ which has a right 
adjoint $R_n$ given by right Kan extension along $H_n$. Since $H_n$ is
fully faithful, so is $R_n$, and so the counit of the adjunction may be
taken to be an identity. For a simplicial object $X$ we write $\Cosk_nX$
for the right Kan extension of its restriction along $H_n$, and
$c:X\to \Cosk_nX$ for the unit map. Then $\Cosk_nX$ is called the 
$n$-coskeleton of $X$, and $X$ is said to be $n$-coskeletal if $h^n$ is 
invertible. We shall be particularly interested in the maps 
$c_n:X_n\to(\Cosk_1X)_n$.

We shall see that the nerve of a category is always 2-coskeletal, while
the 2-nerve of a bicategory is always 3-coskeletal.

\section{Nerves of categories}\label{sect:nerves}

In this section we briefly recall some standard material on nerves of 
categories.

We write $\Delta$ for the category of finite non-empty ordinals and 
order-preserving maps; as usual we write \n for the ordinal
$\{0<1<\ldots<n\}$. Each ordinal can be seen as a category, and 
this provides a fully faithful inclusion functor $J:\Delta\to\Cat_1$
of $\Delta$ in the category $\Cat_1$ of categories and functors.

A simplicial set is a presheaf $X:\Delta\op\to\Set$. We follow the 
usual practice of writing $X_n$ for the image under $X$ of \n.

The singular functor of the inclusion $J:\Delta\to\cat$ is the 
functor $N:\Cat_1\to[\Delta\op,\Set]$ sending a 
category $C$ to its {\em nerve} $NC=\Cat_1(J-,C)$, where 
$\Cat_1(J-,C)$ is the functor sending an ordinal \n to the set
$\Cat_1(J\n,C)$ of all functors from \n to $C$; in other words,
to the set $C_n$ of composable $n$-tuples of morphisms in $C$ (with 
a ``$0$-tuple'' understood to mean just an object). In particular
there are maps $d_0,d_1:C_1\to C_0$ sending a morphism to,
respectively, its codomain and its domain. The maps $d_0,d_2:C_2\to C_1$
are the projections, while $d_1:C_2\to C_1$ is given by composition.

Given any simplicial object $X:\Delta\op\to\E$ in a category \E with finite 
limits, we can form the pulback $X_1\t_{X_0}X_1$
of $d_0,d_1:X_1\to X_0$, and then since $d_0d_2=d_1d_0$ there is a 
map $S_2:X_2\to X_1\t_{X_0}X_1$ induced by $d_0$ and $d_2$. The reason
for the letter ``S'' is that this map, and the $S_n$ described below
have often been called ``Segal maps'', since this approach to coherence
goes back to \cite{Segal74}. Of course there is some 
ambiguity in the notation $X_1\t_{X_0}X_1$, since the maps $X_1\to X_0$
involved are not recorded. We exacerbate this, by abbreviating 
$X_1\t_{X_0}X_1$ to $X^2_1$ (this will {\em always} denote the pullback
constructed in this way). 

A special feature of the simplical sets which are nerves of categories
is that the map $S_2:X_2\to X^2_1$ is invertible; this is precisely the
fact that $X_2$ is the set of composable pairs. Similarly, for an 
arbitrary simplicial object $X$, we write $X^n_1$ for the limit of the 
diagram 
$$
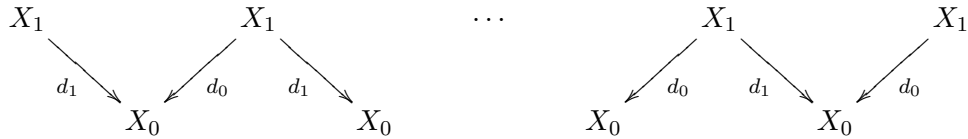
\xymatrix{
X_1 \ar[dr]_{d_1} && X_1 \ar[dl]^{d_0}\ar[dr]_{d_1} && \ldots &&
X_1 \ar[dr]_{d_1}\ar[dl]^{d_0} && X_1 \ar[dl]^{d_0} \\
& X_0 && X_0 && X_0 && X_0 }$$
in which there are $n$ copies of $X_1$, and $S_n:X_n\to X^n_1$ for 
the evident induced map. Once again, for the nerve of a category,
this $S_n$ is invertible.

It is a straightforward but important calculation that the functor
$N:\Cat_1\to[\Delta\op,\Set]$ is fully faithful. (For example, the
fact that a functor preserves composition is encoded in the fact of 
naturality with respect to the map $\delta_1:[1]\to[2]$ whose image
in a simplicial set is $d_1:X_2\to X_1$.) The fact that the nerve
functor is fully faithful may alternatively be expressed by saying that
$J:\Delta\to\Cat_1$ is {\em dense}; but it is probably easiest to 
check the fully faithfulness directly.

There are also smaller subcategories of \cat which are dense. Let
\Deltac denote the subcategory $\Delta$ generated by the objects
$[0]$, $[1]$, and $[2]$, and the morphisms $\delta_0,\delta_1:[0]\to[1]$,
$\sigma_0:[1]\to[0]$, and $\delta_0,\delta_1,\delta_2:[1]\to[2]$ (but not
the degeneracy maps $[2]\to[1]$. Write \Jc for the (non-full)
inclusion $\Deltac\to\cat$. Then $\Jc$ induces a functor 
$\Nc:\cat\to[\Deltacop,\Set]$ sending a category $C$ to the restriction
of its nerve $NC$ to \Deltacop, and once again this functor \Nc is
fully faithful. It follows, by Proposition~\ref{prop:density}, that the 
nerve of a category $C$ is the right Kan
extension along the inclusion $\Deltacop\to\Delta\op$ of $\Nc C$.

Similarly, if \Deltatwo is the full subcategory of $\Delta$ containing 
the objects $[0]$, $[1]$, and $[2]$, then the functor 
$\Ntwo:\cat\to[\Deltatwoop,\Set]$, given by the nerve followed by restriction 
to \Deltatwoop, is fully faithful. Furthermore, by 
Proposition~\ref{prop:density} once again, $N$ itself can be recovered as 
the composite
$$\xymatrix{
\cat \ar[r]^-{\Ntwo}& [\Deltatwoop,\Set] \ar[r]^{\Ran_H} & [\Delta\op,\Set] }$$
where $\Ran_H$ is given by right Kan extension along the inclusion
$H:\Deltatwoop\to\Delta\op$. Thus the nerve of a category is 
always 2-coskeletal (see Section~\ref{sect:singular}).


A simplicial set is the nerve of a category if and only if 
$S_n:X_n\to X^n_1$ is invertible for all $n>1$. (Under the
reasonable definition of $S_1$ and $S_0$, these last are always invertible.)

Since $\Cat_1$ is cocomplete, it follows for general reasons 
(see Section~\ref{sect:singular}) that
the nerve functor has a left adjoint, which sends a simplicial set
$X$ to the weighted colimit $X*J$. Explicitly, the objects of $X*J$
are the elements of $X_0$, while the morphisms are generated by
the elements of $X_1$ subject to relations encoded in $X_2$; the higher
simplices are not needed to calculate $X*J$, essentially because nerves
of categories are 2-coskeletal.


\section{The 2-nerve construction}\label{sect:2-nerves}

We now turn to the case of bicategories. Every category may be seen
as a locally discrete bicategory (that is, a bicategory in which the
only 2-cells are identities). As observed in the introduction, if \bicatl
denotes the category of bicategories and normal lax functors, then 
the (fully faithful) inclusion $H:\Delta\to\bicatl$ induces a fully faithful
map $\bicatl(H,1):\bicatl\to\SSet$ (and so is dense). This construction
might be called the 1-nerve of the bicategory.

Instead, we shall describe a 2-nerve construction. This requires a
2-category \NHom of bicategories. An object of the 2-category will be a 
bicategory, and a morphism will be a normal homomorphism. Given normal 
homomorphisms $F,G:\A\to\B$, there can be a 2-cell from $F$ to $G$
only if $F$ and $G$ agree on objects; a 2-cell then consists of a 2-cell
$\alpha f:Ff\to Gf$ in \B for every 1-cell $f:A\to B$ in \A, subject to the 
following three conditions. First of all, the $\alpha f$ must be natural
in $f$, in the sense that if $\rho:f\to g$ is a 2-cell in \A, then 
$\alpha g.F\rho = G\rho.\alpha f$. Secondly $\alpha f$ must be an 
identity 2-cell if $f$ is an identity 1-cell. Thirdly, if $f:A\to B$ and 
$g:B\to C$ constitute a composable pair in \A, then the diagram
$$\xymatrix{
Fg.Ff \ar[r]^{\alpha g.\alpha f} \ar[d]_{\phi_{f,g}} & 
Gg.Gf \ar[d]^{\psi_{f,g}} \\
F(gf) \ar[r]_{\alpha(gf)} & G(gf) }$$
of 2-cells in \B must commute, where $\phi$ and $\psi$ are the 
pseudofunctoriality
isomorphisms for $F$ and $G$. Such a 2-cell is called an {\em icon}, since
it is precisely an Identity Component Oplax Natural transformation from $F$
to $G$ --- the $\alpha$ are the 2-cells expressing the oplax naturality of
the identity maps $FA\to GA$. 

\begin{remark}
In the important special case where the
bicategories \A and \B have only one object, so that they may be regarded
as monoidal categories, with $F$ and $G$ then becoming strong monoidal
functors, an icon is precisely a monoidal natural transformation.
\end{remark}

\begin{remark}
As pointed out to us by Bob Par\'e, the 2-category \NHom
can be seen as living within the 2-category \LxDbl of pseudo double
categories, lax double functors, and horizontal transformations,
studied by Grandis and Par\'e in \cite{GrandisPare05}. From this point of view,
it is the restriction from pseudo double categories to bicategories
(seen as pseudo double catetgories in which all horizontal arrows
are identities) that leads to the restriction to transformations
whose components are identities.  See the last paragraph of 
\cite[Section~2.2]{GrandisPare05}.
\end{remark}

Every category can be seen as a locally discrete bicategory --- that is, a
bicategory with no non-identity 2-cells. Seen in this way, \cat becomes a
full sub-2-category of \NHom --- a normal homomorphism between locally
discrete bicategories is just a functor between the corresponding categories,
and there are no non-identity icons between such normal homomorphisms.

Thus we can in turn regard $\Delta$ as a full sub-2-category of \NHom, once
again there are no non-identity 2-cells. It is this fully faithful inclusion 
$J:\Delta\to\NHom$ whose singular 2-functor $N=\NHom(J,1):\NHom\to\SCat$
gives our 2-nerve construction.

\begin{remark}\label{rmk:Tamsamani}
This is not the same as the construction described by Tamsamani on 
page~54 of \cite{Tamsamani}: his construction has nothing corresonding
to the coherence condition for normal homomorphisms $\n\to\B$. 
\end{remark}

One of the main results of the paper will be Theorem~\ref{thm:ff} below,
which states that the 2-nerve 2-functor $N:\NHom\to[\Delta\op,\Cat]$ 
is fully faithful.

Just as for categories, it is not necessary to use all of $\Delta$. Write 
\Deltab for the sub-category of $\Delta$ generated by the objects $[0]$, $[1]$,
and $[2]$, and all morphisms between them, as well as the object $[3]$ and the 
four maps $\delta_i:[2]\to[3]$. In fact we shall see that the inclusion 
$H:\Deltab\to\NHom$ is also dense, so that the induced 
$N_b=\NHom(H,1):\NHom\to[\Deltabop,\Cat]$ is fully faithful; and by
Proposition~\ref{prop:density} this will imply that 
$N:\NHom\to[\Delta\op,\Cat]$
is fully faithful; it will likewise imply that the corresponding 2-functor
$\NHom\to[\Deltathreeop,\Cat]$ is fully faithful, where now $\Deltathreeop$
is the full subcategory of $\Delta\op$ consisting of all objects \n with
$n\le 3$.

As a first step to proving that $N_b$ is fully faithful, we describe a little 
more explicitly the 2-nerve 2-functor. We write $\B_n$ for the category 
$\NHom(\n,\B)$ of $n$-simplices of the 2-nerve $N\B$ of \B.

For an ordinal \n and a bicategory \B, a normal homomorhism 
$\n\to\B$ consists of the following data in \B
\begin{itemize}
\item an object $B_i$ for each $i\in\n$
\item a morphism $b_{ij}:B_i\to B_j$ for each $i,j\in\n$ with $i<j$
\item an invertible 2-cell $\beta_{ijk}:b_{jk}b_{ij}\cong b_{ik}$ for
each $i,j,k\in\n$ with $i<j<k$
\end{itemize}
subject to the condition that the diagram
$$\xymatrix @R 1pc {
b_{kl}(b_{jk}b_{ij}) \ar[r]^{b_{kl}\beta_{ijk}} \ar[dd]_{\alpha} &
b_{kl}b_{ik} \ar[dr]^{\beta{ikl}} \\
&& b_{il} \\
(b_{kl}b_{jk})b_{ij} \ar[r]_{\beta_{jkl}b_{ij}} & b_{jl}b_{ij} 
\ar[ur]_{\beta_{ijl}} \\
}$$
commutes for all $i,j,k,l\in\n$ with $i<j<k<l$.

Given another such normal homomorphism $(C,c,\gamma):\n\to\B$, an icon $(B,b,\beta)\to(C,c,\gamma)$ consists of
\begin{itemize}
\item satisfaction of the equation $B_i=C_i$ for each $i\in\n$
\item a 2-cell $\phi_{ij}:b_{ij}\to c_{ij}$ for each $i,j\in\n$, $i<j$
\end{itemize} 
such that the diagram 
$$\xymatrix{
b_{jk}b_{ij} \ar[r]^{\beta_{ijk}} \ar[d]_{\phi_{jk}\phi_{ij}} &
b_{ik} \ar[d]^{\phi_{ik}} \\
c_{jk}c_{ij} \ar[r]_{\gamma_{ijk}} & c_{ik} 
}$$
commutes for all $i,j,k\in\n$ with $i<j<k$. 

Suppose now that $X:\Deltabop\to\Cat$ has $X_0$ discrete and that \B is 
a bicategory. We consider what it is to give a morphism $F:X\to\Nb\B$ in
$[\Deltabop,\Cat]$. If $x$ and $y$ are in $X_0$, we write $X(x,y)$ for
the fibre over $(x,y)$ of the map $X_1\to X_0\t X_0$ induced by $d_1$ and 
$d_0$.

The category $\B_0$ is discrete; its objects are the objects of \B and it has 
no non-identity morphisms, thus $F_0$ simply assigns to each $x\in X_0$
an object $Fx$ of \B.

An object of $\B_1$ is a morphism of \B, while a morphism of $\B_1$ is a 2-cell
in \B. Thus $\B_1$ is the coproduct of the hom-categories $\B(A,B)$ as $A$ and 
$B$ range over all the objects of \B. The face maps $d_0,d_1:\B_1\to\B_0$ give 
the codomain and the domain objects of a 1-cell or 2-cell. The degeneracy map
$s_0:\B_0\to\B_1$ sends an object to the identity 1-cell on the object. Thus to
give $F_1:X_1\to\B_1$, compatible with the face maps $d_0$ and $d_1$, is to 
give a functor $F:X(x,y)\to\B(Fx,Fy)$ for all $x,y\in X_0$.

An object of $\B_2$ consists of morphisms $b_{01}:B_0\to B_1$, 
$b_{12}:B_1\to B_2$, and $b_{02}:B_0\to B_2$, and an invertible 2-cell 
$\beta:b_{12}b_{01}\to b_{02}$. A morphism of $\B_2$ consists of three
2-cells of \B, satisfying the coherence condition given
above. The face maps $d_0,d_1,d_2:\B_2\to\B_1$ pick out the three sides of
the 2-simplex ($d_i$ picks out the map $b_{jk}$, where $i$, $j$, and
$k$ are distinct). The degeneracy maps $s_0,s_1:\B_1\to\B_2$ send a 1-cell 
$b:B\to B'$ to the 2-simplices
defined using the identity isomorphisms $1_{B'}f\cong f$ and $f 1_B\cong f$.
Thus to give $F_2:X_2\to\B_2$ compatible with the face maps is to give, for each 
object $\xi$ of $X_2$, an invertible 2-cell $F_2\xi:Fd_0\xi.Fd_2\xi\cong Fd_1\xi$ in \B,
as in 
\begin{equation}\label{F2}\xymatrix{
& Fd_0d_2\xi \ar[ddr]^{Fd_0\xi} \\
{}\rrtwocell<\omit>{~~F_2\xi} && {} \\
Fd_1d_2\xi \ar[uur]^{Fd_2\xi} \ar[rr]_{Fd_1\xi} && Fd_0d_0\xi }
\end{equation}
natural with respect to the 1-cells in $X_2$. Compatibility with respect to the degeneracy
maps asserts that $F_2s_0f$ and $F_2s_1f$ are the 2-simplices arising from the 
identity isomorphisms $1.Ff\cong Ff\cong Ff.1$ in \B.

Finally to give $F_3:X_3\to\B_3$, compatible with the degeneracy maps, is to assert
that for each object $\Xi$ of $X_3$, the diagram 
\begin{equation}\label{F3}
\xymatrix  @R1pc @C4pc {
Fx_{23}(Fx_{12}.Fx_{01}) \ar[r]^{Fx_{23}.F_2\xi_{012}} \ar[dd]_{\alpha}&
Fx_{23}.Fx_{02} \ar[dr]^{F_2\xi_{023}} \\
&& Fx_{03} \\
(Fx_{23}.Fx_{12})Fx_{01} \ar[r]_{F_2\xi_{123}.Fx_{01}} & Fx_{13}.Fx_{01} 
\ar[ur]_{F_2\xi_{013}} 
}\end{equation}
in \B commutes, where $x_{ij}$ is $d_k d_l\Xi$ for a suitable choice
of $k$ and $l$, and $\xi_{ijk}$ is $d_l \Xi$ for a suitable $l$, while $\alpha$
is the associativity isomorphism.

We record this as:

\begin{proposition}
If $X:\Deltabop\to\Cat$ has $X_0$ discrete, and \B is a bicategory, then to 
give a morphism $F:X\to\Nb\B$ in $[\Deltabop,\Cat]$ is to give (i) an object 
$Fx$ of \B for each $x\in X_0$, (ii) a functor $F:X(x,y)\to\B(Fx,Fy)$ for each 
$x,y\in X_0$, with $Fs_0x=1_{Fx}$, (iii) an invertible 2-cell 
$F_2\xi:Fd_0\xi.Fd_2\xi\cong Fd_1\xi$, as in \eqref{F2} above, for each object 
$\xi\in X_2$, such that (iv)  $F_2\xi$ is natural in $\xi$,
(v) the $F_2s_i\xi$ are the identity isomorphisms, and (vi) \eqref{F3} holds 
for all $\Xi\in X_3$.
\end{proposition}

Suppose now that $F,G:X\to\Nb\B$ are two morphisms in $[\Deltabop,\Cat]$. What 
is it to give a 2-cell (modification) between them? A similar analysis to that 
above gives:

\begin{proposition}
If $X:\Deltabop\to\Cat$ with $X_0$ discrete, \B is a bicategory, and 
$F,G:X\to\Nb\B$, then to give a 2-cell $F\to G$ is (i) to assert that $Fx=Gx$ 
for all $x\in X_0$, (ii) to give a 2-cell $\phi f:Ff\to Gf$ in \B, for every 
$f\in X_1$, such that (iii) $\phi f$ is natural in $f$, (iv) $\phi s_0 x$ 
is an identity 2-cell for each $x\in X_0$, and the diagram 
$$\xymatrix{
Fd_0\xi.Fd_2\xi \ar[rr]^{\phi d_0\xi.\phi d_2\xi} \ar[d]_{F_2\xi} && 
Gd_0\xi.Gd_2\xi \ar[d]^{G_2\xi} \\
Fd_1\xi \ar[rr]_{\phi d_1\xi} && Gd_1\xi}$$
of 2-cells in \B commutes for all objects $\xi$ of $X_2$. 
\end{proposition}

We now specialize the last two propositions to the case where $X$ too has the 
form $\Nb\A$ for a bicategory \A. In the resulting description of a morphism 
$\Nb\A\to\Nb\B$, we see that (i) amounts to the assignment of an object $FA$ 
of \B for each object $A$ of \A, and (ii) amounts to a functor 
$F:\A(A,B)\to\B(FA,FB)$ for all objects $A$ and $B$, with $F1_A=1_{FA}$ for
all $A$. If $f:A\to B$ and 
$g:B\to C$ are 1-cells in \A, then the identity 1-cell on $gf$ determines a 
2-simplex $\xi\in\Nb\A_2$ with faces $f$, $g$, and $gf$, and then $F_2\xi$, 
for this $\xi$, is an invertible 2-cell $\phi_{f,g}:Fg.Ff\cong F(gf)$. 
Essentially by the Yoneda lemma, to give the $F_2\xi$ as in (iii) satisfying 
naturality as in (iv) is just to give such $\phi_{f,g}$, natural in $f$ and 
$g$; then a general 2-simplex $(f,g,\xi:gf\to h)$ must be sent to the 2-simplex
$(Ff,Fg,F_2\xi:Fg.Ff\to Fh)$, where $F_2\xi$ is now the composite
$$\xymatrix{Fg.Ff \ar[r]^{\phi_{f,g}} & F(gf) \ar[r]^{F\xi} & Fh.}$$
Finally (v) and (vi) assert precisely that these $\phi$ satisfy the 
normalization and 3-cocycle conditions to make $F$ into a normal homomorphisms 
from \A to \B. Similarly, if $G:\Nb\A\to\Nb\B$ is another morphism in 
$[\Deltabop,\Cat]$, then to give a 2-cell $F\to G$ is precisely to give an 
icon between the corresponding normal homomorphisms. This proves:

\begin{theorem}
The 2-functor $\Nb:\NHom\to[\Deltabop,\Cat]$ is fully faithful, or equivalently
the inclusion $\Deltab\to\NHom$ is dense.
\end{theorem}

As an immediate consequence we have, by Proposition~\ref{prop:density}:

\begin{theorem}\label{thm:ff}
The 2-nerve 2-functor $N:\NHom\to\SCat$ is fully faithful, or
equivalently the inclusion $\Delta\to\NHom$ is dense.
\end{theorem}

\section{The 2-dimensional universal algebra point of view}\label{sect:2DUA}

As we said in the introduction, the following two sections are not 
needed in the rest of the paper, and can be omitted on a first
reading.

A \Cat-graph \cite{Wolff} has objects $X,Y,Z,\ldots$, with ``hom-categories'' 
$\G(X,Y)$ for each pair of objects $X$ and $Y$. With the obvious notion of 
morphism, this defines a category, which is in fact locally finitely 
presentable. Here, however, we want to make it into a 2-category. Given 
\Cat-graph morphisms $M,N:\G\to\H$, a 2-cell $M\to N$ exists only if $M$ and 
$N$ agree on objects, in which case it consists of a natural transformation 
$$\xymatrix @C 5pc {
~~~~~\G(X,Y) \rtwocell^M_N & \H(MX,MY)}$$
for each pair of objects $X$ and $Y$. These objects, morphisms, and 2-cells now
form a 2-category \CG, which is itself locally finitely presentable, in the 
sense of \cite{Kelly-amiens}.

There is an evident forgetful 2-functor $U:\Hom\to\CG$, and it is a routine 
exercise to give a presentation, in the sense of \cite{KP}, for a finitary 
2-monad $T$ on \CG for which \Hom is the 2-category \talg of (strict) 
$T$-algebras, (pseudo) $T$-morphisms, and $T$-transformations. To see this,
let \bd{n} denote the \Cat-graph with objects $0,1,\ldots,n$, with 
$\bd{n}(i,j)=1$ if $i<j$ and all other hom-categories empty, and let 
$i:\bd0\to\bd{n}$ denote
the map sending $0$ to $i$. Finally let $I$ denote the 
\Cat-graph with  objects 0 and 1, with $I(0,1)$ the 
``free-living isomorphism'', and all other hom-categories empty; thus
there are two isomorphic maps from $0$ to $1$. Then to 
make a \Cat-graph into a bicategory, one must equip it with operations
$$\xymatrix @R 0.5pc {
\CG(\bd2,\G) \ar[r]^{M} & \CG(\bd1,\G) \\
\CG(\bd0,\G) \ar[r]^{j} & \CG(\bd1,\G) \\
\CG(\bd3,\G) \ar[r]^{\alpha} & \CG(I,\G) \\
\CG(\bd1,\G) \ar[r]^{\lambda} & \CG(I,\G) \\
\CG(\bd1,\G) \ar[r]^{\rho} & \CG(I,\G)
}$$
specifying composition, identities, associativity isomorphisms, and left and
right identity isomorphisms, subject to equations between derived operations,
which specify such things as the domain and codomain of composites, and the
coherence condition for the associativity isomorphism. For example the
domains and codomains of composites are specified by commutativity of
the diagrams 
$$\xymatrix{
\CG(\bd2,\G) \ar[r]^{M} \ar[dr]_{\CG(0,\G)} & 
\CG(\bd1,\G) \ar[d]^{\CG(0,\G)} &
\CG(\bd2,\G) \ar[r]^{M} \ar[dr]_{\CG(2,\G)} &
\CG(\bd1,\G) \ar[d]^{\CG(1,\G)} \\
& \CG(\bd0,\G) && \CG(\bd0,\G) }$$
while the domain of the associativity isomorphism is specified by
commutativity of 
$$\xymatrix{
\CG(\bd3,\G) \ar[r]^{\alpha} \ar[d]_{M_1} & 
\CG(I,\G) \ar[d]^{\CG(d,\G)} \\
\CG(\bd2,\G) \ar[r]_{M} & \CG(\bd1,\G) }$$
wherein $d:\bd1\to I$ is one of the identity-on-object inclusions, and
$M_1:\CG(\bd3,\G)\to\CG(\bd2,\G)$ is the map representing ``composing
the first two maps of a composable triple''. The latter is uniquely determined
by commutativity of
$$\xymatrix{
\CG(\bd3,\G) \ar[r]^{M_1} \ar[d]_{\CG(p',\G)} & 
\CG(\bd2,\G) \ar[d]^{\CG(p,\G)} &
\CG(\bd3,\G) \ar[r]^{M_1} \ar[dr]_{\CG(r,\G)} &
\CG(\bd2,\G) \ar[d]^{\CG(q,\G)} \\
\CG(\bd2,\G) \ar[r]_{M} & \CG(\bd1,\G) && \CG(\bd1,\G) }$$
wherein $p:\bd1\to\bd2$ and $p':\bd2\to\bd3$ are inclusions, 
$q:\bd1\to\bd2$ sends an object $i\in\bd1$ to $i+1\in\bd2$, and 
$r:\bd1\to\bd3$ sends $i$ to $i+2$.

Although our main 
interest is in (normal) homomorphisms, and so as usual it is the pseudo
morphisms of $T$-algebras which are most important, the strict morphisms
of $T$-algebras are also of considerable theoretical importance, and 
they are precisely the strict homomorphisms of bicategories.

As a consequence of the fact that \Hom has the form \talg, we may 
deduce, thanks to \cite{BKP}:

\begin{theorem}\label{thm:talg}
The 2-category \Hom has products, inserters, and equifiers, and therefore has 
all bicategorical limits. It also has bicategorical colimits. A homomorphism 
of bicategories is an equivalence if and only if the underlying morphism of 
\Cat-graphs is an equivalence. An icon is invertible if and only if the 
underlying 2-cell in \CG is invertible. If \B is a bicateory, and 
$M:\G\to U\B$ an equivalence in \CG, we may ``transport'' the bicategory 
structure to obtain a bicategory \A and an equivalence $F:\A\to\B$ with 
$U\A=\G$ and $UF=M$. If $G:\C\to\B$ is a homomorphism, and $\rho:N\to UG$ is 
an invertible 2-cell in \CG, we may ``transport'' the homomorphism structure 
to obtain a homomorphism $H:\C\to\B$ and an invertible icon $\sigma:H\to G$ 
with $UH=N$ and $U\sigma=\rho$.
\end{theorem}
 
\section{Some other 2-categories of bicategories}

The 2-category \NHom has a full sub-2-category \Catnps consisting of the 
2-categories. (Here ``nps'' is short for ``normal pseudofunctor'': this is
the name often given to normal homomorphisms between 2-categories.) 
Thus the inclusion $\Catnps\to\NHom$ is fully faithful; we shall see 
that it is also biessentially surjective on objects, and so a biequivalence. 
There is also a larger 2-category \Hom, whose objects are the bicategories, 
but with arbitrary homomorphisms (not necessarily normal) as 1-cells. It is 
fairly straightforward to extend the definition of icons, so as to allow icons 
between arbitrary homomorphisms; the only slight subtlety is that rather than 
asking $\alpha 1_A$  be an identity 2-cell when $1_A$ is an identity 1-cell, 
one rather asks for $\alpha 1_A$ to be suitably compatible with the identity 
constraints $F1_A\cong 1_{FA}$ and $G1_A\cong 1_{GA}$. The inclusion of \NHom 
in \Hom is bijective on objects and locally fully faithful; we shall see
that it is also locally an equivalence, and so a biequivalence. Finally there 
is a full sub-2-category \Catps of \Hom consisting of the 2-categories, and 
once again the inclusions $\Catps\to\Hom$ and $\Catnps\to\Catps$ are biequivalences.

To see these facts, we first describe the equivalences and the invertible 
2-cells in \Hom.

\begin{lemma}
Let $F,G:\A\to\B$ be homomorphisms of bicategories. An  icon $\alpha:F\to G$ 
is invertible if and only if each $\alpha f:Ff\to Gf$ is invertible (in other 
words, if the oplax natural transformation is actually pseudonatural). A 
2-cell in \NHom, \Catps, or \Catnps is an isomorphism if and only if it is one 
in \Hom.
\end{lemma}

\proof
By Theorem~\ref{thm:talg}, the icon $\alpha$ is invertible if and only if the 
underlying 2-cell $U\alpha$ in \CG is invertible, but clearly this says 
precisely that each $\alpha f$ is invertible. The results for \NHom, \Catps, 
and \Catnps are immediate since these are all locally full sub-2-categories of 
\Hom.
\endproof

As a consequence we have:

\begin{proposition}
Every homomorphism $F:\A\to\B$ of bicategories is isomorphic in \Hom to a 
normal homomorphism.
\end{proposition}

\proof
Consider the underlying morphism $UF:U\A\to U\B$ in \CG. There is an evident 
morphism $M:U\A\to U\B$ defined like $UF$ except on the identity 1-cells, 
which are sent to the corresponding identity 1-cells in \B. There is an 
invertible 2-cell $\rho:M\to UF$ in \CG which is the identity except on 
identity 1-cells, where it is the canonical isomorphism $F1_A\cong 1_{FA}$. We 
may transport the homomorphism structure to obtain a homomorphism $G:\A\to\B$ 
and an isomorphism $\sigma:G\to F$ with $UG=M$ and $U\sigma=\rho$. Clearly $G$ 
is in fact a normal homomorphism.
\endproof

\begin{lemma}
Let \A and \B be bicategories. A homomorphism $F:\A\to\B$ is an equivalence in 
\Hom if and only if it is bijective on objects and induces an equivalence of 
hom-categories.  A 1-cell in \NHom, \Catps, or \Catnps is an equivalence if and
only if it is one in \Hom.
\end{lemma}

\proof
The statement about \Hom is more or less immediate from the fact that $F$ is 
an equivalence if and only if the underlying morphism $UF$ of \Cat-graphs is 
an equivalence. The main point of interest is that if $G$ is equivalence 
inverse to $F$, then $GF\cong 1$ and $FG\cong 1$ in \Hom, which forces $GF$ 
and $FG$ to act as the identity on objects, and so for $F$ to be bijective on 
objects. The case of  \Catps follows immediately from that of \Hom; while those
of \NHom and \Catnps now follow using the proposition.
\endproof

It is well-known that every bicategory is biequivalent to a 2-category, and 
that this biequivalence may be chosen to be bijective on objects
\cite{MacLane-Pare}. But a 
biequivalence which is bijective on objects is precisely an equivalence in 
\Hom. This now proves:

\begin{theorem}
Each of the inclusions $\NHom\to\Hom$, $\Catps\to\Hom$, and $\Catnps\to\Hom$ is
a biequivalence of 2-categories.
\end{theorem}
 
The 2-monad of the previous section can be modified so that \talg
becomes not \Hom but \Catps. On the other hand, we can modify the base 
2-category \CG by asking each hom-category $\G(X,X)$ to have a chosen object, 
preserved by the morphisms and 2-cells, and the resulting 2-category \RCG is 
once again locally finitely presentable. There are suitable finitary 2-monads 
on \RCG for which the resulting 2-categories \talg are respectively \NHom and 
\Catnps.

We leave to the reader the resulting modifications of Theorem~\ref{thm:talg} 
dealing with \NHom, \Catps, and \Catnps. As a further consequence of
\cite{BKP} we have:


\begin{theorem}\label{thm:biadjoint}
The singular functor $N=\NHom(J,1):\NHom\to\SCat$ has a left biadjoint.
\end{theorem}

\proof
Let $T$ be the finitary 2-monad on \RCG for which \talg is \NHom. Just
as in the case of \Hom, a {\em strict} morphism of $T$-algebras is precisely
a strict homomorphism of bicategories. We write \NHoms or \talgs for the 
sub-2-category consisting of the strict homomorphisms. The inclusion
$I:\NHoms\to\NHom$ has a left adjoint $L$; this follows once again from 
\cite{BKP}, or can be proved directly.
Now the composite $NI:\NHoms\to\SCat$ can be written as $\NHom(J,I)$, which
by the adjunction $L\dashv I$ is just $\NHoms(LJ,1)$; that is, the singular
functor of $LJ:\Delta\to\NHoms$. Now \NHoms is cocomplete {\em as a 
2-category}, and so $NI=\NHoms(LJ,1)$ has a 2-adjoint $F$; thus by
\cite[Theorem~5.1]{BKP}, the composite $IF$ is biadjoint to to $N$. 

The left biadjoint $IF$ is in fact a 2-functor, and the unit $1\to NIF$ is 
2-natural; neither of these facts is true for a general biadjunction.
\endproof

We shall also see below that if we restrict the codomain of $N$ we can 
obtain a very special left 2-adjoint, which is almost a 2-equivalence.

\section{Properties of 2-nerves}\label{sect:properties}

We have already seen that the 2-nerve $N\B$ of a bicategory is the right
Kan extension of a 2-functor $\Nb:\Deltabop\to\Cat$ along the inclusion
$\Deltabop\to\Delta\op$, and so in particular that it is 3-coskeletal.
We have also seen that $N\B_0$ is discrete.

For each $n>1$, the Segal map $S_n:\B_n\to\B^n_1$ is a surjective 
equivalence. When $n=2$, for example, this says (i) for a composable
pair $f:A\to B$ and $g:B\to C$, there exist a morphsim $h:A\to C$ and
an invertible 2-cell $\phi:gf\cong h$, and (ii) given data as above,
and also $f':A\to B$, $g':B\to C$, $h':A\to C$, and $\phi':g'f'\cong h'$,
then for each pair $\alpha:f\to f'$ and $\beta:g\to g'$ of 2-cells, there is
a unique 2-cell $\gamma:h\to h'$ for which the evident pasting diagram
(involved in the definition of morphisms of 2-simplices) commutes.
Of course for (i), we may take $h=gf$ and $\phi$ to be the identity,
while for (ii), we may (and must!) take $\gamma$ to be the composite
$\phi'.\beta\alpha.\phi^{-1}$.

The fact that these Segal maps are surjective equivalences will be 
important in the following section, where we turn to the notions of
weak 2-category due to Tamsamani and to Simpson.

We say that a functor $p:E\to B$ is a {\em discrete isofibration} or
{\em dif}, if for each $e\in E$ and each isomorphism $\beta:b\to pe$ in $B$,
there exists a unique isomorphism $\epsilon:e'\to e$ in $E$ with 
$p\epsilon=\beta$ (and so also $pe'=b$). (This property might equally
be called ``unique transport of structure'' or ``unique invertible-path
lifting''.)

Recall that we write $c_n$ for the $n$-component of the canonical map
$X\to\Cosk_{n-1}X$ from a simplicial object $X$ to its $n-1$-coskeleton.
Our next observation is that $c_2$ is a dif for $X=N\B$. A 2-simplex
in the 1-coskeleton of $N\B$ consists of three maps
$f:A\to B$, $g:B\to C$, and $h:A\to C$ in \B; thus the fact that $c_2$ is
a dif amounts to the (evident) fact that if $\phi:gf\cong h$, and we are given 
isomorphisms $\alpha:f'\to f$, $\beta:g'\to g$, and $\gamma:h'\to h$,
then there is a unique way to paste these together to obtain an isomorphism
$\phi':g'f'\to h'$.

Similarly $c_3$ is a dif: this amounts to the (equally
evident) fact that for an
$n$-simplex consisting of morphisms $x_{ij}:X_i\to X_j$ and invertible
2-cells $\xi_{ijk}:x_{jk}x_{ij}\to x_{ik}$, given invertible 2-cells
$\zeta_{ij}:x'_{ij}\to x_{ij}$ for each $i<j$, when one constructs the 
unique induced $\xi'_{ijk}$ guaranteed by the fact that $c_2$ is a
dif, these $\xi'_{ijk}$ fit together to form an $3$-simplex.


We shall see below that a simplicial object $X:\Delta\op\to\Cat$ is 
the 2-nerve of a bicategory if and only if (i) $X_0$ is discrete, 
(ii) $S_n:X_n\to X^n_1$ is an equivalence for all $n>1$, 
(iii) $c_2$ and $c_3$ are difs, and (iv) $X$ is 3-coskeletal.

In the remainder of this section we establish two results which 
will be used in the comparison between bicategories and Tamsamani
weak 2-categories to which we turn in the following section. 
Up until now, the only sort of morphisms between functors
$\Delta\op\to\Cat$ have been the (2-)natural transformations. But since
\Cat is a 2-category, it is possible, and indeed reasonable, to consider
also pseudonatural transformations. There is a 2-category 
$\Ps(\Delta\op,\Cat)$ of (2-)functors, pseudonatural transformations, and
modifications, but we shall also be interested in the sub 
2-category $\NPs(\Delta\op,\Cat)$ containing only those pseudonatural
transformations $X\to Y$ for which the pseudonaturality isomorphism with
respect to each epimorphism in $\Delta$ is an identity. We call such
a pseudonatural transformation a {\em normal pseudonatural}, since the
strict naturality with respect to epimorphisms is closely related to 
normality of pseudofunctors/homomorphisms. Notice that ``normal pseudonatural 
transformation'' is not ascribed any meaning in general, only in this special 
case of pseudonatural transformations between functors from $\Delta\op$ to 
\Cat.

The first of these two results says that we can ``normalize'' pseudonatural
transformations whose domain is a 2-nerve.

\begin{proposition}\label{prop:normalization}
Any morphism $f:N\A\to X$ in $\Ps(\Delta\op,\Cat)$, with $N\A$ the
2-nerve of a bicategory \A, is isomorphic to a normal pseudonatural
transformation.
\end{proposition}

\proof
We shall define inductively functors $g_n:(N\A)_n\to X_n$ equipped with 
natural isomorphisms $\psi_n:f_n\cong g_n$ such that the composite
$$\xymatrix{
g_{n+1}s_i \ar[r]^{\psi^{-1}_{n+1}s_i} & f_{n+1}s_i \ar[r]^{f_{s_i}} & 
s_i f_{n} \ar[r]^{\psi_n s_i} & s_i g_n }$$
is an identity for each $i$; here the isomorphism 
$f_{s_i}:f_{n+1}s_i\cong s_i f_n$ is the pseudonaturality isomorphism.
Then the pseudonaturality isomorphisms for $f$ can be transported across
the isomorphisms $\psi_n$ to give pseudonaturality isomorphisms for
the $g_n$, and the composite displayed above will be precisely the
induced pseudonaturality isomorphism $g_{n+1}s_i\cong s_i g_n$; thus
$g$ will become normal pseudonatural, in such a way that $\psi$ is
an invertible modification.

We take $g_0$ to be $f_0$, and $\psi_0$ to be the identity. To define
$g_{n+1}$ and $\psi_{n+1}$, it suffices to choose, for each object
$x\in(N\A)_{n+1}$, an object $g_{n+1}x\in X_{n+1}$ and an isomorphism
$\psi_{n+1}x:f_{n+1}x\cong g_{n+1}x$, such that the composite
$$\xymatrix{
g_{n+1}s_i y \ar[r]^{\psi^{-1}_{n+1}s_i y} & f_{n+1}s_i y \ar[r]^{f_{s_i y}} & 
s_i f_{n}y \ar[r]^{\psi_n s_i y} & s_i g_n y}$$
is an identity for each $i$; then $g_{n+1}$ becomes a functor and $\psi_{n+1}$
a natural isomorphism.

If $x$ is non-degenerate, we take $g_{n+1}x$ to be $f_{n+1}x$ and 
$\psi_{n+1}x$ to be the identity. If $x=s_i y$, we take 
$g_{n+1}x=s_i g_n y$, and $\psi_{n+1}x$ to be the composite 
$\psi_n s_i y.f_{s_i y}:f_{n+1}s_i y\to s_i g_n y = g_{n+1}s_i y$.
The only thing to check is that this is well-defined.
Now $s_i$ is a section, so there can be at most one $y$ with $s_i y=x$;
but it is possible that $x=s_i y=s_j z$, with $j<i$. It is at this
stage that we use the fact that the domain of $f$ is the 2-nerve of
a bicategory. For in this case we necessarily have $y=s_j w$ and $z=s_{i-1}w$
for some $w\in(N\A)_{n-1}$: this boils down to 
the fact that in a bicategory the left and right identity isomorphisms
$1_A 1_A\cong 1_A$ must agree. Now 
$s_i g_n y = s_i g_n s_j w = s_i s_j g_{n-1} w = s_j s_{i-1} g_{n-1} w =
s_j g_n s_{i-1} w = s_j g_n z$, and so $g_{n+1}x$ is indeed well-defined,
and the well-definedness of $\psi_{n+1}x$ is similar.
\endproof

\begin{remark}
It is clear from the proof that one could relax the assumption that 
the domain is the 2-nerve of a bicategory. What is really used is 
that certain commutative squares of degeneracy maps are actually
pullbacks.
\end{remark}

The second result asserts a kind of ``fibrancy'' or 
``coflexibility'' property of simplicial objects $X:\Delta\op\to\Cat$
of the form $N\B$ for a bicategory \B. (See \cite{qm2cat} for more
about the relationship between flexibility and cofibrancy.) It
shows that if we have a {\em normal} pseudonatural transformation
whose {\em codomain} is the 2-nerve of a bicategory, then we can 
replace it by an isomorphic 2-natural transformation. Combined with
the previous result, this will imply that any pseudonatural transformation
between 2-nerves of bicategories is isomorphic to a 2-natural transformation.

It was shown in \cite{BKP} that the inclusion 
$[\A,\B]\to\Ps(\A,\B)$ admits a left adjoint whenever \A is a small 
2-category and \B a cocomplete one. For suitable choices of \A and \B
it follows that $[\Delta\op,\Cat]\to\Ps(\Delta\op,\Cat)$ admits both
a left and a right adjoint. A straightforward variant of this (which
is still in fact a special case of the main theorem of \cite{BKP})
shows that likewise the inclusion $[\Delta\op,\Cat]\to\NPs(\Delta\op,\Cat)$ 
admits both adjoints. It is the right adjoint of the latter inclusion which
concerns us here; it sends a functor $X:\Delta\op\to\Cat$ to
another such functor $X^+$; and by the universal property of the adjunction
combined with the Yoneda lemma we see that
\begin{align*}
X^+_n &\cong [\Delta\op,\Cat](\Delta(-,n),X^+) \\
     &\cong \NPs(\Delta\op,\Cat)(\Delta(-,n),X)
\end{align*}
and now the component at $X$ of the unit of the adjunction is the 
map $j:X\to X^+$ whose $n$-component is the inclusion of 
$[\Delta\op,\Cat](\Delta(-,n),X)$ in $\NPs(\Delta\op,\Cat)(\Delta(-,n),X)$.

The {\em counit} is a normal pseudonatural $p:X^+\to X$ with $pj=1$,
and the arguments of \cite{BKP} show that $jp\cong 1$, and so that
$j$ is an equivalence in $\NPs(\Delta\op,\Cat)$. We shall say that $X$ is 
{\em coflexible} (but {\em fibrant} would also be a good name) if
$j:X\to X^+$ has a retraction $r$ in $[\Delta\op,\Cat]$; it then follows,
as in \cite{BKP}, that $r\cong rjp=p$, and so $jr\cong jp\cong 1$, and so 
finally that $j$ is an equivalence in $[\Delta\op,\Cat]$.

Alternatively, rather than relying on \cite{BKP}, one could {\em define}
$X^+$ by $X^+_n=\NPs(\Delta\op,\Cat)(\Delta(-,n),X)$ and $j$ as the 
inclusion, and then prove directly that it has the properties stated above.

\begin{theorem}\label{thm:coflexible}
The 2-nerve $N\B$ of any bicategory \B is coflexible.
\end{theorem}

\proof
First we make slightly more explicit the description of $X^+$, 
when $X=N\B$.

A normal pseudonatural from $\Delta(-,n)$ to $X$ consists of an
object $\xi\in X_n$ equipped with, for each non-identity monomorphism 
$\delta_i:\m\to\n$
in $\Delta$, an object $\xi_\delta$ and an isomorphism 
$u_{\delta_i}:d_i\xi\cong\xi_\delta$ in $X_m$. A morphism between
two such objects is just a morphism between their underlying objects
in $X_n$.

In particular, since $[0]$ has no subobjects, $X^+_0$ is just $X_0$.
On the other hand, while $[1]$ does have two subobjects (the two
maps $\delta_0,\delta_1:[0]\to [1]$), the category $X_0$ has no non-identity
isomorphisms, and so once again $X^+_1$ is just $X_1$.

When it comes to $X^+_2$ things become more interesting. Once again,
there are no non-identity isomorphisms in $X_0$, but there are 
non-identity isomorphisms in $X_1$, and three non-identity monomorphisms
$[1]\to[2]$. Thus an object of $X^+_2$ consists of an object $\xi\in X_2$,
equipped with an isomorphism $u_i:d_i\xi\cong \xi_i$ in $X_1$, for 
$i=0,1,2$. A morphism in $X^+_2$ is just a morphism between the underlying
objects in $X_2$. The inclusion $j:X_2\to X^+_2$ equips $\xi\in X_2$
with identity isomorphisms. The degeneracies $X^+_1\to X^+_2$ are
induced by $j$ from the degeneracies $X_1\to X_2$.

An object in $X^+_3$ consists of an object $\Xi\in X_3$, equipped with
isomorphisms $v_i:d_i\Xi\cong\Xi_i$ in $X_2$ for $i=0,1,2,3$; as well
as isomorphisms $w_d:d\Xi\cong\Xi_d$ in $X_1$, where $d$ runs though
each of the six monomorphisms $[1]\to[3]$ in $\Delta$. The face map
$d_i:X^+_3\to X^+_2$ sends such an object $(\Xi,(v_i),(w_d))$ to
the 2-simplex in $X^+$ consisting of $\Xi_i$, equipped with 
the isomorphism $x_{ij}:d_j\Xi_i\cong\Xi_{ij}$ given by 
$$\xymatrix{
d_j\Xi_i \ar[r]^-{d_jv^{-1}_i} & d_jd_i\Xi = d\Xi \ar[r]^-{w_d} & \Xi_d}$$
where $d=d_jd_i$; we shall sometimes write $\Xi_d=\Xi_{ij}$.

We shall now construct the desired retraction $r:X^+\to X$, using the
fact that $X=N\B$, and the description of morphisms into such objects
of $[\Delta\op,\Cat]$.

For $n=0$ and $n=1$, the map $j:X_n\to X^+_n$ is the identity, and so
we take $r$ to be the identity; clearly $r:X^+_1\to X_1$ is
compatible with degeneracies. We define $r_2:X^+_2\to X_2$ on objects
to take $(\xi,(u_i))$ to the 2-simplex in $X_2$ given by 
$$\xymatrix{
\xi_2.\xi_0 \ar[r]^{u^{-1}_2.u^{-1}_0} & d_2\xi.d_0\xi \ar[r]^{\xi} &
d_1\xi \ar[r]^{u_1} & \xi_1.}$$
This is compatible with the face maps by construction, and compatible
with the degeneracies by definition of degeneracies in $X^+_2$. The
functoriality of $r_2$ is clear.

It remains to check the compatibility condition codified by a 3-simplex.
But this asserts the commutativity in the bicategory \B
of diagram~(\ref{F3}) in Section~\ref{sect:2-nerves}, which asserts
the equality, for each object of $X^+_3$, of 
a parallel pair of arrows in $X_1$. But the construction of this parallel
pair is natural in the objects of $X^+_3$, and holds for objects in the
image of $j:X_3\to X^+_3$, so holds for all objects, since $j$ is an
equivalence.
\endproof

\section{The Tamsamani and Simpson notions of weak 2-category}
\label{sect:Tamsamani}


Let \Tam denote the full sub-2-category of \SCat consisting of those
$X$ for which $X_0$ is discrete and each $S_n:X_n\to X^n_1$ is an equivalence, 
and \Simp the smaller full sub-2-category of those $X$ for which 
moreover each $S_n$ is surjective. We speak of {\em Tamsamani 2-categories}
and {\em Simpson 2-categories}, but in fact Tamsamani used the name 
{\em 2-nerve}, while Simpson used the name {\em easy 2-category}.

We saw in the previous section that if $X$ is the 2-nerve of a bicategory,
then each $S_n:X_n\to X^n_1$ is a surjective equivalence and $X_0$ is discrete.
Thus the 2-nerve of a bicategory is a Simpson 2-category, and so necessarily
a Tamsamani 2-category.

The left biadjoint of Theorem~\ref{thm:biadjoint} also gives biadjoints to the
inclusions of \NHom in each of \Simp and \Tam, but in fact there exist 
2-adjoints. For a Tamsamani 2-category $X$, there is a 
bicategory $GX$, defined in \cite{Tamsamani}, whose objects are the elements 
of $X_0$, and whose 1-cells 
and 2-cells are the objects and morphisms of $X_1$, with vertical composition 
of 2-cells given by the composition law in $X_1$. Since 
$S_2:X_2\to X_1\t_{X_0}X_1$ is an equivalence, we can choose a functor 
$M:X_1\t_{X_0}X_1\to X_1$ and an isomorphism $\sigma:d_1\cong MS_2$,
as in 
$$\xymatrix @R 1pc @C 1pc {
X_2 \ar[rr]^{S_2} \ar[dd]_{d_1} \drtwocell<\omit>{\sigma} && 
X^2_1 \ar[ddll]^{M} \\
& {} \\
X_1 }$$
and this $M$ gives the 
composition of 1-cells and the horizontal composition of 2-cells.

Composing $\sigma$ with the degeneracy maps $s_0,s_1:X_1\to X_2$ gives
$$\xymatrix @R 1pc {
1=d_1s_0 \ar[r]^-{\sigma s_0} & MS_2s_0 = M\binom{1}{s_0d_1} \\
1=d_1s_1 \ar[r]^-{\sigma s_1} & MS_2s_1 = M\binom{s_0d_0}{1}
}$$
which are the identity isomorphisms for the bicategory $GX$. For the 
associativity isomorphism
$M(M\t 1)\cong M(1\t M)$, consider the pasting diagrams
$$
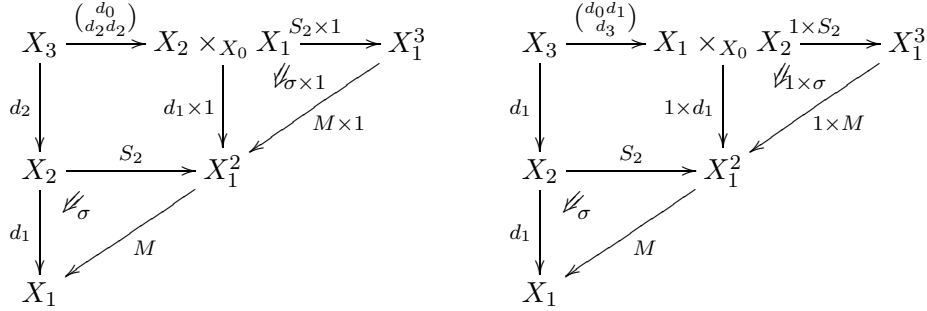
\xymatrix @R 1pc @C 1pc {
X_3 \ar[rr]^-{\binom{d_0}{d_2d_2}} \ar[dd]_{d_2} && 
X_2\t_{X_0}X_1 \ar[rr]^{S_2\t 1} \ar[dd]_{d_1\t1} 
\drtwocell<\omit>{~~\sigma\t1}&& X^3_1 \ar[ddll]^{M\t 1} &&
X_3 \ar[rr]^-{\binom{d_0d_1}{d_3}} \ar[dd]_{d_1} && 
X_1\t_{X_0}X_2 \ar[rr]^{1\t S_2} \ar[dd]_{1\t d_1} 
\drtwocell<\omit>{~~1\t\sigma} && X^3_1 \ar[ddll]^{1\t M}\\
&& & {}& &&
&&  & {} \\
X_2 \ar[rr]^{S_2} \ar[dd]_{d_1} \drtwocell<\omit>{\sigma} && 
X^2_1 \ar[ddll]^{M} && &&
X_2 \ar[rr]^{S_2} \ar[dd]_{d_1} \drtwocell<\omit>{\sigma} && 
X^2_1 \ar[ddll]^{M} \\
& {} & && &&
 & {}\\
X_1 && && && X_1}$$
The left hand composite of the two diagrams are equal, and the top composite
of each diagram is just the equivalence $S_3:X_3\to X^3_1$; thus there is
a unique invertible $\alpha:M(M\t 1)\cong M(1\t M)$ which when pasted onto
the left diagram gives the right diagram, and we take this $\alpha$ to be 
the associativity isomorphism. The coherence condition for associativity
is checked using the fact that $S_4$ is an equivalence: the two pasting
composites which are to be proven equal are each pasted onto the right hand
side of the following diagram (in which the operation ``$~\t_{X_0}~$'' has
been written as juxtaposition)
$$
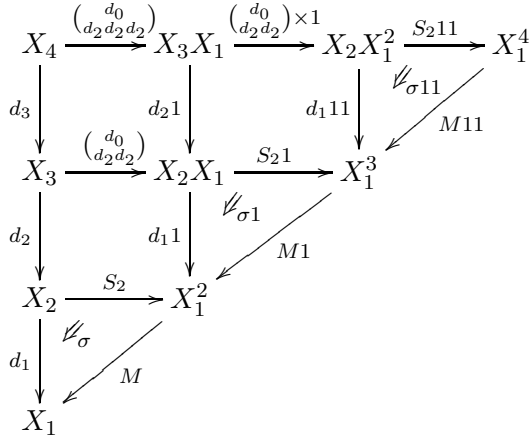
\xymatrix @R 1pc @C 1pc {
X_4 \ar[rr]^{\binom{d_0}{d_2d_2d_2}} \ar[dd]_{d_3} && 
X_3X_1 \ar[rr]^{\binom{d_0}{d_2d_2}\t1} \ar[dd]_{d_2 1}&&
X_2X^2_1 \ar[rr]^{S_211} \ar[dd]_{d_111} \drtwocell<\omit>{~~\sigma11} 
&& X^4_1 \ar[ddll]^{M11} \\
&&&&& {} \\
X_3 \ar[rr]^{\binom{d_0}{d_2d_2}} \ar[dd]_{d_2} && 
X_2X_1 \ar[rr]^{S_2 1} \ar[dd]_{d_1 1} \drtwocell<\omit>{~\sigma1} 
&& X^3_1 \ar[ddll]^{M1} \\
&&& {} \\
X_2 \ar[rr]^{S_2} \ar[dd]_{d_1} \drtwocell<\omit>{\sigma} 
&& X^2_1 \ar[ddll]^{M} \\
& {} \\
X_1
}$$
and in each case, after three steps, one gets the same result, namely
$$\xymatrix @R 1pc {
X_4 \ar[rr]^{\binom{d_0d_0d_0}{d_2}} \ar[dd]_{d_1} && 
X_1X_3 \ar[rr]^{1\binom{d_0d_1}{d_3}} \ar[dd]_{1d_1} &&
X^2_1 X_2 \ar[rr]^{11S_2} \ar[dd]_{11d_1} \drtwocell<\omit>{~~11\sigma} 
&& X^4_1 \ar[ddll]^{11M} \\
&&&&& {} \\
X_3 \ar[rr]^{\binom{d_0d_1}{d_3}} \ar[dd]_{d_1} && 
X_1X_2 \ar[rr]^{1 S_2} \ar[dd]_{1 d_1} \drtwocell<\omit>{~1\sigma} 
&& X^3_1 \ar[ddll]^{1M} \\
&&& {} \\
X_2 \ar[rr]^{S_2} \ar[dd]_{d_1} \drtwocell<\omit>{\sigma} 
&& X^2_1 \ar[ddll]^{M} \\
& {} \\
X_1
}$$
Now both these pasting composites are invertible, and the arrow across
the top is the equivalence $S_4:X_4\to X^4_1$, so the desired result 
follows. A similar but easier argument establishes the coherence
for the identities.


Having constructed the bicategory $GX$, we can now of course take its
2-nerve $NGX$, which we temporarily name $X'$. We are going to construct
a morphism $u:X\to NGX=X'$. Clearly $X'_0=X_0$ and
$X'_1=X_1$, and so we may take $u_0$ and $u_1$ to be the identities.
Now $X'_2$ may be constructed as the {\em pseudolimit}
of the map $M:X^2_1\to X_1$; that is, the universal category equipped
with functors $S'_2:X'_2\to X^2_1$ and $d'_1:X'_2\to X_1$,
and an isomorphism $\rho:MS'_2\cong d'_1$. The components $d'_0$ and $d'_2$
of $S'_2$, along with $d'_1$, are the face maps $X'_2\to X'_1=X_1$. Thus
a map into $X'_2$ is determined by its composite with the face maps
$X'_2\to X_1$ and with $\rho$. By a general property of such 
pseudolimits $S'_2$ must be a surjective equivalence.

We construct $u_2:X_2\to X'_2$ as the unique map compatible with the
face maps $X'_2\to X'_1$, and with $\rho u_2=\sigma$. Notice that
$S'_2u_2=\binom{d'_0}{d'_2}u_2=\binom{d_0}{d_2}=S_2$, and both $S'_2$ and 
$S_2$ are equivalences, thus also $u_2$ is an equivalence.
The degeneracy map $s'_0:X'_1=X_1\to X'_2$ is the unique map satisfying
$d'_0s'_0=d'_1s'_0=1$, $d'_2s'_0=s_0d'_1$, and $\rho s'_0=\sigma s_0$;
similarly $s'_1:X_1\to X'_2$ is the unique map satisfying 
$d'_0s'_1=s'_0d'_0$, $d'_1s'_1=d'_2s'_1=1$, and $\rho s'_1=\sigma s_1$.
It is now straightforward to verify that $u_2$ is compatible with 
the degeneracy maps.

The category $X'_3$ is the pseudolimit of the diagram
$$\xymatrix{
X^3_1 \ar[r]^{M\t 1} \ar[d]_{1\t M} \drtwocell<\omit>{\alpha} & 
X^2_1 \ar[d]^{M} \\
X^2_1 \ar[r]_{M} & X_1 }$$
in other words, the universal category $X'_3$ equipped with 
morphisms $S'_3:X'_3\to X^3_1$ and $K_1,K_2:X'_3\to X^2_1$ and $L:X'_3\to X_1$,
and with invertible 2-cells $\kappa_1:K_1\cong (1\t M)S'_3$, 
$\kappa_2:(M\t 1)S'_3\cong K_2$, and $\lambda_1:L\cong MK_1$. (This may appear 
asymmetric but it is not; there is a uniquely induced isomorphism
$\lambda_2:L\cong MK_2$ suitably compatible with the associativity
isomorphism $\alpha$.) Once again, it is an immediate consequence that
$S'_3$ is a surjective equivalence. 

There is a unique map $u_3:X_3\to X'_3$ with $S'_3u_3=S_3$, $K_1u_3=S_2d_1$,
$K_2u_3=S_2d_2$, $Lu_3=d_1d_2$, $\kappa_1u_3=d_0d_0\t \sigma d_3$,
$\kappa_2u_3=\sigma d_0\t d_2d_2$, while $\lambda_1u_3=\sigma d_1$. Since
$S'_3u_3=S_3$ and $S'_3$ and $S_3$ are equivalences, $u_3$ is an 
equivalence too.

Next we describe the face maps $X'_3\to X'_2$. First of all $d'_1$ is
the unique map induced by $K_1$, $L$, and $\lambda_1$, while similarly
$d'_2$ is that induced by $K_2$, $L$, and $\lambda_2$. The other two 
are only slightly more complicated: $S'_3$, $K_1$, and $\kappa_1$ are 
composed with the projection $X^2_1\to X_1$ onto the first factor,
to obtain the data inducing $d'_3$, while $S'_3$, $K_2$, and $\kappa_2$
are composed with the projection $X^2_1\to X_1$ onto the second factor
to obtain the data inducing $d'_0$. One now verifies that $u_3$ is 
compatible with these face maps: for example, $d'_1u_3$ is the map 
$X_3\to X'_2$ induced by $K_1u_3=S_2d_3$, $Lu_3=d_1d_3$, and
$\lambda_1u_3=\sigma d_1$, but these are precisely the data that
describe $u_2d_1$. The other face maps are treated similarly, and so
we have now defined a map $u$ between the restrictions of $X$ and $X'$
to $\Deltabop$, but since $X'$ is the right Kan extension of its 
restriction to $\Deltabop$, it follows that $u$ extends uniquely to
give a map $X\to X'$.

Since $X$ is a Tamsamani 2-category and $X'$ is the 2-nerve
of a bicategory, the Segal maps $S_n:X_n\to X^n_1$ and $S'_n:X'_n\to X^n_1$
are equivalences; on the other hand, $u_n:X_n\to X'_n$ is clearly 
compatible with the Segal maps, and so it follows that each $u_n$
is an equivalence. A morphism $f:X\to Y$ of simplicial object is
said to be a {\em pointwise equivalence} if, as here, each $f_n$
is an equivalence.

Thus for every Tamsamani 2-category $X$, we have constructed
a pointwise equivalence $u:X\to NGX$ to the 2-nerve of a bicategory $GX$.

Suppose now that $X$ is 3-coskeletal and $c_2$ and $c_3$ are difs.
We shall show that $u_2$ and $u_3$ are isomorphisms, so that $u$ 
is an isomorphism, and $X$ is the 2-nerve of a bicategory (namely $GX$).

We already know that $u_2$ is an equivalence; we must show that it
is bijective on objects. Suppose then that $\phi:gf\cong h$ is an
object of $X'_2$. Since $S_2:X_2\to X^2_1$ is an equivalence, there
exists a $\xi\in X_2$ with isomorphisms $\alpha:f\cong d_2\xi$ and
$\beta:g\cong d_0\xi$. The isomorphism $\sigma:MS_2\xi\cong d_1\xi$ 
involved in the definition of composition in $GX$ goes from
$d_0\xi.d_2\xi$ to $d_1\xi$. Combining it with $\alpha$, $\beta$, and
$\phi^{-1}$ provides an isomorphism $\gamma:h\to d_1\xi$ as in
$$\xymatrix{
h\ar[r]^{\phi^{-1}} & gf \ar[r]^{\beta.\alpha} & d_0\xi.d_2\xi \ar[r]^{\sigma}&
d_1\xi,}$$
and now $\alpha$, $\beta$, and $\gamma$ together constitute an isomorphism 
$\phi':(f,g,h)\cong c_2\xi$ in $(\Cosk_1X)_2$, and so since $c_2$ is a dif, 
there is a unique $\phi'':\xi'\cong\xi$ in $X_2$ with $c_2\xi''=(f,g,h)$ and 
$c_2\phi''=\phi'$; in other words, with $u_2\xi''=(f,g,\phi:gf\cong h)$.
This proves that $u_2$ is bijective on objects, and so an isomorphism.

Once again, we already know that $u_3$ is an equivalence, and must show
that it is bijective on objects. Injectivity is easy: if $\Xi$ and $\Xi'$
are objects of $X_3$ with $u_3\Xi=u_3\Xi'$, then there is a unique 
isomorphism $\Xi\cong\Xi'$ sent by $u_3$ to the identity; but it then 
easily follows that it is sent by $c_3:X_3\to(\Cosk_1X)_3$ to the identity,
and so that it must itself be an identity. As for surjectivity, suppose
now that $\Xi'$ is an object of $X'_3$. We know that there is an isomorphism
$\phi:\Xi'\cong u_3\Xi$ for some $\Xi\in X_3$. Applying the map
$c'_3:X'_3\to(\Cosk_1X')_3$ gives an isomorphism
$c'_3\phi:c'_3\Xi'\cong c'_3u_3\Xi$ in $\Cosk_1X')_3$. But $X'$ and $X$
have the same 1-coskeleton, so $(\Cosk_1X')_3=(\Cosk_1X)_3$ and $c'_3u_3=c_3$;
and now the previous isomorphism may be seen as an isomorphism 
$c'_3\phi:c'_3\Xi'\cong c_3\Xi$ in $(\Cosk_1X)_3$. Since $c_3$ is a 
discrete isofibraion, there is a unique isomorphism 
$\phi_1:\Xi_1\cong\Xi$ in $X_3$ with $c_3\phi_1=c'_3\phi$ (and so 
also $c_3\Xi_1=c'_3\Xi'$). But now $u_3\phi_1:u_3\Xi_1\cong u_3\Xi$ is
an isomorphism in $X'_3$ with the property that 
$c'_3u_3\phi_1=c_3\phi_1=c'_3\phi$, and so since $c'_3$ is also a 
discrete fibration, $u_3\phi_1=\phi$ and $u_3\Xi_1=\Xi'$.
This proves:



\begin{theorem}\label{thm:characterization}
A 2-functor $X:\Delta\op\to\Cat$ is the 2-nerve of a bicategory if 
and only if (i) it is 3-coskeletal, (ii) $X_0$ is discrete, (iii) the Segal
maps $S_n:X_n\to X^n_1$ are equivalences, and
(iv) $c_2:X_2\to(\Cosk_1X)_2$ and $c_3:X_3\to(\Cosk_1X)_3$ are discrete
isofibrations.
\end{theorem}



We now show that $u:X\to NGX$ is the unit of a 2-adjunction
between \NHom and \Tam. To do this, let \B be a bicategory and 
$F:X\to N\B$ a morphism in \Tam (equivalently, in $[\Delta\op,\Cat]$).
We must show that there is a unique morphism $F':NGX\to N\B$ with $F'u=F$;
the two dimensional aspect of the universal property will then follow, 
since \NHom has cotensors, preserved by $N:\NHom\to\Tam$. Now to give 
$F':NGX\to N\B$ is equivalently to give
a normal homomorphism $F'':GX\to\B$. Since $u_0$ and $u_1$ are the identities,
the action of $F''$ on objects, 1-cells, and 2-cells is already determined;
it will remain only to give the pseudofunctoriality isomorphisms.

To give a map $F:X\to N\B$ is to give a function $F=F_0:X_0\to\B_0$,
functors $F:X(x,y)\to\B(Fx,Fy)$, and an invertible 2-cell
$\chi_{\xi}:Fd_2\xi.Fd_0\xi\to Fd_1\xi$ for each object $\xi\in X_2$
such that (i) $\chi_{\xi}$ is natural in $\xi$, (ii) $\chi_{\xi}$ is 
the relevant identity isomorphism whenever $\xi$ is degenerate, and 
(iii) if $\Xi$ is an object of $X_3$, then the diagram
$$\xymatrix{
(Fd_0d_0\Xi.Fd_2d_0\Xi).Fd_2d_3\Xi \ar[r]^-{\chi_{d_0\Xi}.1} \ar@{=}[d] &
Fd_1d_0\Xi.Fd_2d_3\Xi \ar@{=}[r] & 
Fd_0d_2\Xi.Fd_2d_2\Xi \ar[r]^-{\chi_{d_2\Xi}} & Fd_1d_2 \Xi \ar@{=}[dd] \\
(Fd_0d_1\Xi.Fd_0d_3\Xi).Fd_2d_3\Xi \ar[d]_{\alpha} \\
Fd_0d_1\Xi.(Fd_0d_3\Xi.Fd_2d_3\Xi) \ar[r]_-{1.\chi_{d_3\Xi}} & 
Fd_0d_1\Xi.Fd_1d_3\Xi \ar@{=}[r] & 
Fd_0d_1\Xi.Fd_2d_1\Xi \ar[r]_-{\chi_{d_1\Xi}} & Fd_1d_1\Xi }$$
in \B commutes.

How can we extend this to a morphism $NGX\to N\B$? On 0-simplices and
1-simplices there is no change: we still use the same function $F:X_0\to\B_0$
and the same functors $F:X(x,y)\to\B(Fx,Fy)$. When it comes to 2-simplices,
we know that every object $(f,g,\phi:gf\cong h)$ of $NGX$ is isomorphic to
an object of the form $u_2\xi$ for a $\xi\in X_2$, so must be sent to a
2-simplex in $N\B$ isomorphic to the image of $\xi$ under $F$. But 
compatibility with the face maps tells where the faces $f$, $g$, and $h$
must go --- namely to $Ff$, $Fg$, and $Fh$, and now everything else is
uniquely determined. Explicitly, fix $\xi\in X_2$, and an isomorphism
$u_2\xi\cong(f,g,\phi)$, given by $\alpha:d_2\xi\cong f$, 
$\beta:d_0\xi\cong g$, and $\gamma:d_1\xi\cong h$. Then $(f,g,\phi)$
must be sent to the 2-simplex of $N\B$ made up of $Ff$, $Fg$, and
$$\xymatrix @C 3pc {
Fg.Ff \ar[r]^-{F\beta.F\alpha} & Fd_0\xi.Fd_2\xi \ar[r]^-{F_2(\xi)} &
Fd_1\xi \ar[r]^-{F\gamma^{-1}} & Fh}$$
where $F_2(\xi)$ is the image of the 2-simplex $\xi\in X_2$ under the
map $F:X\to N\B$. (Notice that the final result does not depend on the 
choice of $\xi$ or the isomorphism $u_2\xi\cong(f,g,\phi)$.

This proves:

\begin{theorem}
The 2-nerve 2-functor $N:\NHom\to\Tam$, seen as landing in the 2-category
\Tam of Tamsamani 2-categories, has a left 2-adjoint given by $G$.
Since $N$ is fully faithful, the counit $GN\to 1$ is invetible. Each 
component $u:X\to NGX$ of the unit is a pointwise equivalence, and
$u_0$ and $u_1$ are identities.
\end{theorem}

If $u$ were in fact an equivalence in \Tam, then $N$ would be fully
faithful and biessentially surjective, and so a biequivalence.
Since each $u_n:X_n\to NGX_n$ is an equivalence, we can choose 
inverse equivalences $v_n:NGX_n\to X_n$, and these will automatically
become the components of a {\em pseudonatural} transformation 
$v:NGX\to X$, but there is no reason in general why they should be
natural, and so there is no reason in general why $u:X\to NGX$ should
be an equivalence in \Tam. One response would be to ``localize'' \NHom
and \Tam by inverting certain morphisms (and throwing away the 2-cells). 

As is always the case with a full reflective subcategory, if one inverts the 
components of the unit --- in this case the $u:X\to NGX$ --- then one 
recovers the subcategory. One could, however, consider inverting
larger classes of maps, for instance the pointwise equivalences,
or more generally the {\em weak equivalences} (or {\em external equivalences}
in \cite{Tamsamani}): a morphism $f:X\to Y$ is a weak equivalence if
and only if $Gf$ is a biequivalence of bicategories. One can show
that inverting these maps in \NHom and \Tam gives equivalent
categories. In fact, in the case where one uses the weak equivalences, the 
resulting categories are also equivalent to
the homotopy categories of the Quillen model categories \twocat and
\Bicats of \cite{qm2cat,qmbicat}; here \twocat is the category of
2-categories and 2-functors, and \Bicats the category of bicategories
and strict homomorphisms.

A less violent approach than inverting these morphisms is to use a 
simplicial localization, as in \cite{Dwyer-Kan1980a}; this time,
using \cite[Corollary~3.6]{Dwyer-Kan1980a} , one obtains weakly
equivalent simplicial categories after localization.

Here we adopt a more precise and explicit approach, in which we expand
our notion of morphism in \Tam to allow not just natural transformations,
but pseudonatural ones. Let \Tamps be the full sub-2-category of 
$\Ps(\Delta\op,\Cat)$ consisting of the Tamsamani 2-categories.

\begin{theorem}\label{thm:main}
The 2-nerve 2-functor $\NHom\to\Tamps$ is a biequivalence of 2-categories.
An object of \Tamnps is in the image of the functor if and only if it 
satisfies the conditions of Theorem~\ref{thm:characterization}; a morphism
between 2-nerves of bicategories
is in the image if and only if it is (not just pseudonatural but 2-)
natural.
\end{theorem}

\proof
The statements about the image have already been proven.

We already know that the 2-nerve 2-functor is locally fully faithful.
To say that it is locally essentially surjective on objects is to 
say that any pseudonatural transformation $f:N\A\to N\B$ is isomorphic 
to a 2-natural transformation. By Proposition~\ref{prop:normalization},
$f$ is isomorphic to a normal pseudonatural transformation $g$. By the
universal property of $p:N\B^+\to N\B$, there is a unique 2-natural 
$h:N\A\to N\B^+$ for which $g=ph$, and now by Theorem~\ref{thm:coflexible}
there is a 2-natural $r$ isomorphic to $p$, and so $g=ph\cong rh$ with $rh$
2-natural. Thus the nerve 2-functor is locally an equiavlence. It remains to
show that it is biessentially surjective on objects; that is, that
every Tamsamani 2-category is equivalent in \Tamps to the 2-nerve
of a bicategory. But if $X$ is a Tamsamani 2-category then we
have the 2-nerve $NGX$, and the pointwise equivalence $u:X\to NGX$, and
every pointwise equivalence is an equivalence in \Tamps.
\endproof

\begin{remark}
One could also consider the sub-2-category \Tamnps of \Tamps containing
only the normal pseudonatural maps, and this would once again be
biequivalent: this time local essential surjectivity
uses only Theorem~\ref{thm:coflexible}, but one now needs to use
Proposition~\ref{prop:normalization} to prove biessential surjectivity.
\end{remark}

Finally we consider what happens in the one-object case. The full 
sub-2-category of \NHom consisting of the one-object bicategories is
precisely the 2-category of monoidal categories, normal strong
monoidal functors, and monoidal natural transformations. On the other hand,
a one-object Tamsamani weak 2-category is what has sometimes been called
a homotopy monoidal category \cite{Leinster-book}. The 2-adjunction
between \NHom and \Tam 
restricts to the one-object case, and so once again, one obtains a weak
equivalence between the simplicial localizations of the category of 
monoidal categories and strong monoidal functors, and the category 
of homotopy monoidal categories and morphisms thereof. Likewise, 
Theorem~\ref{thm:main} restricts to the one-object case.

\bibliographystyle{plain}

\end{document}